\documentclass[12pt]{amsart}
\usepackage{amscd}
\usepackage{amsmath}
\usepackage{amsfonts}
\usepackage{amssymb}
\usepackage{verbatim}
\input prepictex
\input pictexwd
\input postpictex

\newtheorem{theorem}{Theorem}
\theoremstyle{plain}

\newtheorem{fact}{Fact}

\newtheorem{corollary}{Corollary}

\newtheorem{definition}{Definition}

\newtheorem{lemma}{Lemma}

\newtheorem{proposition}{Proposition}

\numberwithin{equation}{section}

\newcommand{\Ien}{i_1\dots i_{n}}
\newcommand{\Jen}{j_1\dots j_n}

\newcommand{\Ken}{k_1\dots k_{n}}
\newcommand{\Ienn}{i'_1\dots i'_{n}}
\newcommand{\Jenn}{j'_1\dots j'_n}
\newcommand{\Mad}{\left\{0,\dots,M\!-\!1\right\}}
\newcommand{\Madn}{\left\{0,\dots,M\!-\!1\right\}^n}
\newcommand{\ien}{\underline{i}_n}
\newcommand{\jen}{\underline{j}_n}

\newcommand{\ken}{\underline{k}_n}
\newcommand{\lem}{\underline{l}_m}

\begin{document}

\parindent0pt

\title[The difference of random Cantor-sets]
{On the size of the algebraic
difference  of two random Cantor sets}

\author{Michel Dekking}
\address{Michel Dekking, Delft Institute of Applied Mathematics,
 Technical University of Delft,  The Netherlands
\newline
\tt{F.M.Dekking@math.tudelft.nl}}

\author{K\'{a}roly Simon}
\address{K\'{a}roly Simon, Institute of Mathematics, Technical
University of Budapest, H-1529 B.O.box 91, Hungary
\tt{simonk@math.bme.hu}}

 \thanks{2000 {\em Mathematics Subject Classification.} Primary
28A80 Secondary 60J80, 60J85
\\ \indent
{\em Key words and phrases.} Random fractals, Mandelbrot percolation,
 difference of Cantor sets, Palis conjecture, multitype branching
 processes in varying environment, superbranching processes.\\
\indent The research of Dekking was partially supported by  the
National Science Foundation of China \#10371043. The research of
Simon was supported by OTKA Foundation \#T42496 and the NWO-OTKA
common project}

\begin{abstract}
In this paper we consider some families of random Cantor sets on
the line and investigate the question whether the condition that
the sum of Hausdorff dimension is larger than one implies the
existence of interior points in the difference set of two
independent copies. We prove that this is the case for the so
called Mandelbrot percolation. On the other hand the same is not
always true if we apply a slightly more general construction of
random Cantor sets. We also present a complete solution for the
deterministic case.
\end{abstract}

\maketitle

\section{Introduction}

\indent{}Algebraic differences of Cantor sets occur naturally in
the context of the dynamical behavior of diffeomorphisms. From
these studies a conjecture by Palis(\cite{P}) originated, relating
the size of the arithmetic difference $F_2-F_1=\{y-x: x\in F_1,
y\in F_2\}$
 to the Hausdorff dimensions of the two Cantor sets $F_1$ and $F_2$: if
\begin{center}
 $\dim_{\rm H}F_1+\dim_{\rm H}F_2>1$ \end{center}
 then  \emph{generically} it should be true that
$$F_2-F_1 \mbox{ contains an interval}.$$
For generic  dynamically generated \emph{non-linear} Cantor sets
this was proved in 2001 by de Moreira and Yoccoz (\cite{MY}). The
problem is open for generic linear Cantor sets. The problem was put
into a probabilistic context by Per Larsson in his thesis \cite{L},
(see also \cite{L-note}). He considers a two parameter family of
random Cantor sets $F_{a,b}$, and obtains that the Palis conjecture
holds for a set of $a$ and $b$ of full Lebesgue measure.
However Larsson's proof contains
errors and significant gaps. In a forthcoming paper the authors of
the present paper will correct these errors and fill the gaps in Larsson's
proof.
Here we will study Palis' conjecture for a natural class of
random Cantor sets considered e.g. in \cite{FRF}, \cite{DG},
\cite{F-proj} and \cite{FG}. A special member of this class was already
 considered  in 1974 by Mandelbrot (\cite{Man}).

\section{Random Cantor sets}\label{21}
 Given are $M\geq 2$ and
the vector $\mathbf{p}:=\left(p_0,\dots,p_{M-1}\right) \in
[0,1]^M$, in general \emph{not} a probability vector and $p_i=0$
or $1$  are also  allowed.

Let $\mathcal{T}$ be the $M$-adic tree. For each $n$ $\mathcal{T}$
has $M^n$ nodes at level $n$, which we denote by strings
$\ien=\Ien$, where $i_k\in\Mad$ for $k=1,\dots,n$. There is one
node at level 0, the root, denoted $\emptyset$. We consider a
probability measure $\mathbb{P}_\mathbf{p}$ on the space of
labeled trees, i.e., each node $\Ien$ obtains a label $X_{\Ien}$
which will be 0 or 1. The probability measure is defined by
requiring that the $X_{\Ien}$ are independent Bernoulli random
variables, with $\mathbb{P}_\mathbf{p}(X_\emptyset=1)=1$, and for
$n\ge 1$ and $\Ien\in\Madn$
$$\mathbb{P}_\mathbf{p}(X_{\Ien}=1)=p_{i_n}.$$
 In particular, when the $X_{\Ien}$ are i.i.d.---i.e. $p_i=p$ for all
$i$---then $\mathbb{P}_\mathbf{p}$ will generate Mandelbrot
percolation.

The randomly labeled  tree generates a random Cantor set in [0,1]
in the following way. Define
\begin{equation*}
I_{\Ien}:=\left[\frac{i_1}{M}+\frac{i_2}{M^2}+\cdots
 +\frac{i_n}{M^n},\frac{i_1}{M}+\frac{i_2}{M^2}+\cdots+\frac{i_n}{M^n} +
 \frac{i_n+1}{M^n}\right].
\end{equation*}
The $n$-th level approximation $F^n$ of the random Cantor set is a
union of such $n$-th level $M$-adic intervals selected by the sets
$S_n$ defined by
$$ S_n=\{\Ien: X_{i_1}=X_{i_1i_2}=\dots=X_{\Ien}=1\}.$$
The random Cantor set $F$ is
$$F=\bigcap_{n=1}^{\infty } F^n=\bigcap_{n=1}^{\infty } \bigcup_{\;\;\Ien \in S_n}I_{\Ien}. $$
 Let $Z_n=\,$Card$(S_n)$ be the number of non-empty intervals $I_{\Ien}$ in $F^n$ and let
$Z_0:=1$.
 Then $\left(Z_n\right)_{n\in \mathbb{N}}$ is a
branching process with offspring distribution the law of $Z_1$.
Namely, let $\xi _{i}^{(n)}$, for $i,n\geq 1$ be i.i.d. random
variables such that $\xi
_{i}^{(n)}\stackrel{\textrm{d}}{=}Z_1$. Then    
\begin{equation*}
Z_{n+1}:=\left\{%
\begin{array}{ll}
    \xi _{1}^{(n+1)}+\cdots +\xi _{Z_n}^{(n+1)},
    & \hbox{if $Z_n>0$;} \\
    0, & \hbox{if $Z_n=0$}. \\
\end{array}%
\right.
\end{equation*}
Note that
$$\mathbb{E}_\mathbf{p}(Z_1)=\mathbb{E}_\mathbf{p}(X_0+\dots+X_{M-1})=p_0+\dots+p_{M-1}.$$
Therefore the branching process will almost surely die out---and
$F$ will be empty---if this expectation is smaller than 1. Hence
we will  assume from now on that
\begin{equation}\label{nonem}
\sum\limits_{k=0}^{M-1}p_k>1.
\end{equation}
The expectation also determines  the Hausdorff dimension
$\dim_{\rm H}F$ of $F$; it is well known (\cite{FRF} or
\cite{MWRF}) that:
\begin{fact}\label{8}
$\dim_{\rm H}F =\log\left(\sum\limits_{k=0}^{M-1}p_k\right)\big / \log M$
almost surely on $F\ne \emptyset$.
\end{fact}

\medskip

\section{Differences of Random Cantor sets}\label{DiffRCS}

Let $F_1,F_2$ be two independent copies of the random Cantor set
$F$ above. From now on $\mathbb{P}$ will denote the product
probability $\mathbb{P}_\mathbf{p}\times\mathbb{P}_\mathbf{p}$.
  Let $F_1^n$ and $F_2^n$ be the corresponding $n$-th level
  approximants of $F_1$ and $F_2$, so
$$
F_i:=\bigcap \limits_{n=1}^{\infty }F_{i}^{n}, \mbox{ for } i=1,2.
$$
Our aim here is to investigate whether the difference set \newline
$$F_2-F_1=\left\{y:\exists\, x_i\in F_i, y=x_2-x_1\right\}$$
contains an interval. It is immediate that

\begin{fact}\label{9}
For a set $A\subset \mathbb{R}^2$ we denote the projection of $A$
on the $y$ axis along lines having a $45^{\circ}$ angle with the $x$
axis by $\mbox{Proj}_{45^{\circ}}(A)$. Then
\begin{equation*}
F_2-F_1=\mbox{Proj}_{45^{\circ}}\left(F_1\times F_2\right).
\end{equation*}
\end{fact}
\noindent  In this way, if $\dim_{\rm H}F<\frac{1}{2}$ then
$\dim_{\rm H}\left(F_2-F_1\right)<1$, so it does not contain any
interval. By Fact \ref{8}, this happens if and only if
$\sum\limits_{k=0}^{M-1}p_k<\sqrt[]{M}$. So, we may hope to find
an interval in $F_2-F_1$ only if the following condition holds:
\begin{equation}\label{11}
\dim_{\rm H}F_1+\dim_{\rm H}F_2>1, \mbox{ that is }
\sum\limits_{k=0}^{M-1}p_k>\sqrt[]{M}.
\end{equation}

Define $p_{M+j}=p_j$ for $j=0,1,\dots,M$. Now we can define the
 \emph{cyclic autocorrelations} $\gamma_k$ by
$$\gamma_k:=\sum\limits_{j=0}^{M-1}p_jp_{j+k} \quad \mbox{for}\quad k=0,\dots ,M.$$


\begin{theorem}\label{12}
Conditional on $F_1,F_2\not=\emptyset $, we have
\begin{description}
    \item[(a)] If $\gamma_k>1$ for all $k$ then $F_2-F_1$ contains an
    interval almost surely.
    \item[(b)] If there exists an $k\in \Mad$ such that  $\gamma_k$ and $\gamma_{k+1}$
    are both less than $1$ then
    $F_2-F_1$ almost surely does not contain any intervals.
\end{description}
\end{theorem}

 In the case of the Mandelbrot percolation  all $p_i=p$  for some
 $0\le p \le 1$. In this case $\gamma_k=Mp^2$ for all $k$. With
 Fact~\ref{8} we obtain the following corollary.

\begin{corollary}\label{Mand}
The Palis conjecture holds for Mandelbrot percolation. That is, if
$F$ is  Mandelbrot percolation, then $\dim_{\rm H}F_1+\dim_{\rm
H}F_2<1$ implies that $F_2-F_1$ almost surely contains no
interval, and  $\dim_{\rm H}F_1+\dim_{\rm H}F_2>1$ implies that
$F_2-F_1$ almost surely does contain an interval (conditional on
$F_1,F_2$ being non-empty).
\end{corollary}

\section{Comments on Theorem 1}

\subsection{Exceptional behaviour} It can happen that Condition (\ref{11}) holds
but almost surely, $F_2-F_1$ does not contain an interval. Let
$M=3$ and for a small number $\varepsilon >0$ (say, $\varepsilon<1/4$)
let $p_0=1,p_1=0,p_2=1-\varepsilon $. (This is
almost the triadic Cantor set with the difference that the second
interval is chosen with probability less than one.) Then Condition
(\ref{11}) holds, but $\gamma_1=\gamma_2=1-\varepsilon <1$, so
almost surely there is no interval in $F_2-F_1$. That is, the
so-called Palis Conjecture does not hold.

\subsection{Scope of the theorem}\label{scope}
In the general case it can happen that for some $k$, $\gamma_k<1$
but $\gamma_{k+1}>1$. In this case our theorem is inconclusive
(see Section~\ref{101rho} for a further discussion). However, if
$M=3$ then $\gamma_0\geq \gamma_1=\gamma_2$. Thus, if
$p_0p_1+p_1p_2+p_2p_0>1$ then $F_2-F_1$ almost surely contains an
interval given that
 $F_1,F_2$ are non empty. On the other hand, if $p_0p_1+p_1p_2+p_2p_0<1$ then
$F_2-F_1$ does not contain any interval almost surely.

\subsection{The deterministic case}
In the case that all $p_i\in\{0,1\}$ we have a complete answer to
the question whether $F_2-F_1$ contains an interval or not. This
will be given in Section~\ref{determ}.

\subsection{A generalisation}

Theorem~\ref{12} remains true when we consider the difference set
of two independent Cantor sets $F_1$ and $F_2$ generated by two
different $p$-vectors of the same length; the autocorrelations
simply have to be replaced by  cross correlations. Assume that the
probabilities for $F_1$ are $p_0,\dots ,p_{M-1}$ and for $F_2$ the
probabilities are $q_0,\dots ,q_{M-1}$. Then to get
$$
\dim_{\rm H}F_1+\dim_{\rm H}F_2>1
$$
we need to assume that
$$
\sum\limits_{i=0}^{M-1}p_i\cdot \sum\limits_{j=0}^{M-1}q_j>M.
$$
The cross correlations are:
$$
\gamma _k:=\sum\limits_{j=0}^{M-1}q_jp_{j+k}
$$
With this all calculations will be the same except for a small
adaptation of the proof of Lemma~\ref{48}. The obvious
generalization of Corollary~\ref{Mand} remains true.

\section{Counting triangles}

Before we start the proof of Theorem \ref{12} we would like to introduce some
notation. Since it is easier to study  $90^\circ$ projections
 we rotate the $[0,1]\!\times [0,1]$ square by $45^\circ$
in the positive direction and translate it, so that its horizontal
diagonal, let us call it $J$, is the
$\left[-\frac12\sqrt{2},\frac12\sqrt{2}\right]$ interval on the
$x$ axis. Let this transformation be called $\varphi$, and let
\begin{equation*}
Q:=\varphi([0,1]\!\times\![0,1]),
\quad\Lambda^n:=\varphi(F_1^n\times F_2^n),
\quad\Lambda:=\varphi(F_1\times F_2).
\end{equation*}
In this way instead of the $45^\circ$ degree projection
$\mbox{Proj}_{45^\circ}$ of $F_1\times F_2$ to the $y$ axis, it is
equivalent to consider the orthogonal projection of $\Lambda $ to
$J$.

\begin{figure}[b!]
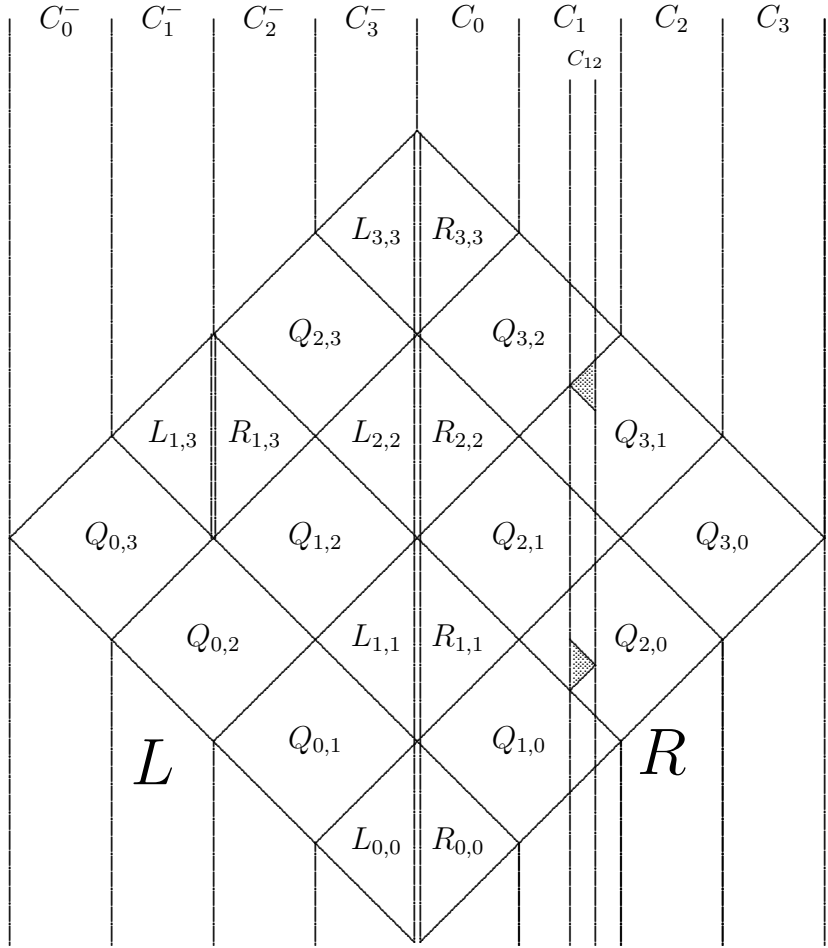

$$
\beginpicture
\setcoordinatesystem units <0.1\textwidth,0.1\textwidth>
\setplotarea x from -1 to 9, y from -4 to 5.6
\setlinear
\linethickness=13pt
\setlinear
\plot  0 0   4 4  /
\plot  0 0   3.97  -3.97 / \plot  4.03  -3.97  8 0  /
\plot  1 -1    5 3 / \plot  2 -2    6 2 / \plot  3 -3    7 1  /
\plot  1.98 2    1.98 0 / \plot  2.02 2    2.02 0 /
\plot  0 -4    0 0 / \plot  1 -4    1 -1 / \plot  2 -4    2 -2 /
\plot  3 -4    3 -3 / \plot  3.97 -3.97    3.97 3.97 / \plot  4.03
-3.97    4.03 3.97 / \plot  5 -4    5 -3 / \plot  6 -4    6 -2 /
\plot  7 -4    7 -1 / \plot  5 -4    5 -3 / \plot  6 -4    6 -2 /
\plot  7 -4    7 -1 / \plot  8 -4    8 5.1 / \plot  8  0    8 5.1
/ \plot  7  1    7 5.1 / \plot  6  2    6 5.1 / \plot  5  3    5
5.1 / \plot  4 4     4 5.1 / \plot  0  0    0 5.1 / \plot  1  1
1 5.1 / \plot  2  2    2 5.1 / \plot  3  3    3 5.1 /
\plot  5.5 4.5    5.5 -4 / \plot  5.75 4.5    5.75 -4 /
\plot 5.5 1.5  5.75 1.25 / \plot 5.5 -1   5.75 -1.25 / \plot 5.5
-1.5 5.75 -1.25 /
\plot  1  1    5 -3 / \plot  2  2    6 -2 / \plot  3  3    7 -1 /
\plot  4  4    8  0 /
\put {\huge $R$} at 6.4 -2.1 \put {\huge $L$} at 1.4 -2.2
\put {$C_{0}$} at 4.5 5.1 \put {$C_{1}$} at 5.5 5.1 \put {$C_{2}$}
at 6.5 5.1 \put {$C_{3}$} at 7.5 5.1 \put {$C_{0}^-$} at 0.5 5.1
\put {$C_{1}^-$} at 1.5 5.1 \put {$C_{2}^-$} at 2.5 5.1 \put
{$C_{3}^-$} at 3.5 5.1
\put {$Q_{0,3}$} at 1 0 \put {$Q_{3,2}$} at 5 2 \put {$Q_{2,1}$}
at 5 0 \put {$Q_{1,0}$} at 5 -2
\put {$L_{1,3}$} at 1.6 1 \put {$R_{1,3}$} at 2.4 1 \put
{$Q_{0,2}$} at 2 -1 \put {$Q_{3,1}$} at 6.2 1 \put {$Q_{2,0}$} at
6.2 -1 \put {$Q_{2,3}$} at 3 2 \put {$Q_{1,2}$} at 3 0 \put
{$Q_{0,1}$} at 3 -2 \put {$Q_{3,0}$} at 7 0
\put {$R_{3,3}$} at 4.4 3 \put {$R_{2,2}$} at 4.4 1 \put
{$R_{1,1}$} at 4.4 -1 \put {$R_{0,0}$} at 4.4 -3 \put {$L_{3,3}$}
at 3.6 3 \put {$L_{2,2}$} at 3.6 1 \put {$L_{1,1}$} at 3.6 -1 \put
{$L_{0,0}$} at 3.6 -3
\put { \tiny $C_{12}$} at 5.6  4.7
\setshadegrid span <0.03cm> 
\vshade 5.5 1.5 1.5 5.75 1.25 1.75 /
\setshadegrid span <0.03cm> 
\vshade 5.5 -1.5 -1.  5.75 -1.25 -1.25 /
\endpicture
$$
\caption{The case $M=4$. The square $Q$ split into two triangles
$L$ and $R$ with level 1 squares, triangles and columns. Also
shown is level 2 column $C_{12}$ with a $\Delta$-pair.}
\label{fig:Delta}
\end{figure}

The image under  $\varphi$ of the  square $I_{\Ien}\!\times
I_{\Jen}$ is denoted $Q_{\Ien,\Jen}$. Every $n$-th level square
$Q_{\Ien,\Jen}$ is divided into two congruent triangles by its
vertical diagonal. The one which is on the left side is denoted
$L_{\ien,\jen}$. We call $L_{\Ien,\Jen}$ an $n$-th level
$L$-triangle. The other part of the square $Q_{\ien,\jen}$ is
denoted $R_{\ien,\jen}$.  In the same way, we divide the  square
$Q$ into two triangles $L$  and $R$ , as in Figure
\ref{fig:Delta}. Note that  $\Lambda$ satisfy a symmetry property:
if we replace  $(x,y)$ by $(-x,y)$ (i.e. $(\ien,\jen)$ by
$(\jen,\ien)$ at level $n$) then $L\cap \Lambda$ is mapped to
$R\cap \Lambda$ and vice versa. Moreover, since this corresponds
to replacing $F_1\times F_2$ by $F_2\times F_1$, $\mathbb{P}$ is
invariant for this mirroring. It follows that properties that we
deduce for $R\cap \Lambda$ will also hold for $L\cap \Lambda$. For
this reason, and to simplify the statements, several of the
following results are formulated for the $R$-triangle only.

%


The orthogonal projection (any projection from now on will be meant
to be orthogonal) of the $n$-th level $L$- and $R$-triangles to
$\left[0,\frac12\sqrt{2}\right]$ are $M^n$ intervals of length
$\frac12\sqrt{2}\cdot M^{-n}$. We denote them in the following
way:%
$$
J_{\Ken}:=\frac12\sqrt{2}\cdot I_{\Ken}.
$$
The intervals in  $\left[-\frac12\sqrt{2},0\right]$ will be denoted
$$
J_{\Ken}^-:=J_{\Ken}-\tfrac12\sqrt{2}.
$$

Now we introduce the appropriate vertical columns intersecting
triangle $L$ respectively $R$:  we write
\begin{equation*}\label{13}
 C_{\Ken}^-:=J_{\Ken}^-\times\mathbb{R}, \qquad    C_{\Ken}:=J_{\Ken}\times\mathbb{R}.
\end{equation*}

We are going to count the number of $L$- and $R$-triangles in
these columns. The idea is that as long as there are $L$- and
$R$-triangles in $C_{\Ken}$, then the interval $J_{\Ken}$ is in
the projection of $\Lambda^n$. Let for  $U,V\in
\left\{L,R\right\}$ the number of level 1 $V$-triangles in
$\Lambda^1\cap C_k^-$ (if $U=L$), respectively $\Lambda^1\cap C_k$
(if $U=R$), generated by a level 0 $U$-triangle be denoted by
$Z^{UV}(k)$. So for instance
\begin{equation*}
Z^{LR}(k)=\#\left\{(i,j): Q_{i,j}\subset \Lambda^1,R_{i,j}\subset  C_k^-\right\}.
\end{equation*}
More generally, we denote by $Z^{UV}(\ken)$ the number of level
$n$ $V$-triangles in $\Lambda^n\cap C_{\ken}^-$ (respectively $\Lambda^n\cap C_{\ken}$)
generated  by a level 0 $U$-triangle. Let
$$
\mathcal{M}(\ken):=\begin{bmatrix}
  \mathbb{E}Z^{LL}(\ken) & \mathbb{E}Z^{LR}(\ken) \\
  \mathbb{E}Z^{RL}(\ken) & \mathbb{E}Z^{RR}(\ken)
\end{bmatrix}.
$$
Then from the definition one can easily check that
\begin{equation}\label{33}
\mathcal{M}(\Ken)=\mathcal{M}(k_1) \cdots \mathcal{M}(k_n).
\end{equation}
In the context of branching processes this is an obvious property:
for a fixed sequence $(k_1,k_2,\dots,k_n,\dots)$ the process
$(Z^{UV}(\ken))$ is a two type branching process in a varying
environment with neighbour interaction. Actually we will be needing another process
that is not a multi type branching process, but a superbranching
process (\cite{DG}) in varying environment, which also has neighbour interaction.
The reason we need this process is that an $R$-triangle will have very few
offspring (next level $L$- and $R$-triangles) in the right most subcolumns
in which it occurs, and possibly a lot in the left most columns. To balance this
asymmetry we will pair $L$-triangles and $R$-triangles.

Let $n\ge 1$. Any  pair
$(L^n,R^n)=(L_{\Ien,\Jen},R_{{\Ienn},\Jenn})$ of \emph{disjoint}
$n$-th level $L$-triangles and $R$-triangles with $L^n\subset
C_{\ken}$ and $R^n\subset C_{\ken}$ for some $\ken$ is called a
level $n$ $\Delta$-pair (cf.~Figure \ref{fig:Delta}). Disjoint
means that they are not allowed to share an edge, so that their
offspring distributions can not interact.
We try to find as many $\Delta$-pairs as possible in a column, where an
 $L$-triangle or $R$-triangle is allowed to belong to at most  one
$\Delta$-pair.
It can be  proved by mathematical induction that for $m\geq 3$
we can find $m$  $\Delta$-pairs as soon as
 $m$ $L$-triangles, and at least $m$ $R$-triangles occur,  or vice versa.
We then say that $m$ $\Delta$-pairs occur in the column.  To analyze the process of
$\Delta$-pairs we consider the number of $V$-triangles generated
in columns $ C_{\ken}^-$  and $ C_{\ken}$ by the level 0 triangles
$L$ \emph{and} $R$. This is equal to
$$Z^V(\ken):=Z^{LV}(\ken)+Z^{RV}(\ken).$$


Note (cf. Figure \ref{fig:Delta}) that for each $k \in \Mad$:
$$L_{i,j}\in C_k \Leftrightarrow i-j=k+1\; (\text{mod} M),\quad
    R_{i,j}\in C_k \Leftrightarrow i-j=k\; (\text{mod} M),$$
and that this also holds for the $C^-_k$ columns.

Since
$$\mathbb{P}(Q_{i,j}\in \Lambda^1)=\mathbb{P}(I_i\in F^1_1,
I_j\in F^1_2)=\mathbb{P}_\mathbf{p}(I_i\in
F^1_1)\mathbb{P}_\mathbf{p} (I_j\in F^1_2)=p_ip_j, $$
 we obtain that for  $k \in \Mad$
\begin{equation}\label{36}
\mathbb{E}Z^L(k)=\gamma_{k+1}, \qquad
\mathbb{E}Z^R(k)=\gamma_{k}.
\end{equation}
 This is the reason that in the statement and the proof of Theorem~\ref{12} the number
$$
\gamma:=\min\limits_{V\in \left\{L,R\right\},0\leq k\leq
M-1}\mathbb{E}Z^V(k)=\min\limits_{0\leq k\leq M-1}\gamma_k
$$
has an important role.

We will need in the proof of Lemma~\ref{pos} that there is a positive probability
that the number of $\Delta$-pairs  grows  exponentially fast in all columns.
When all $p_i$ are positive this is trivial. The following lemma deals with the
case where some $p_i$ may be zero.

\begin{lemma}\label{37}
For any $n$
$$\mathbb{P}\left(Z^L(\ken)\ge\gamma^n\; {\rm and} \; Z^R(\ken)\ge\gamma^n
\;\;{\rm for\; all} \;\ken\right)>0.$$
\end{lemma}
\begin{proof}
 We want to know the maximal number of level $n$ $L$-triangles and
 $R$-triangles
  that occur in the columns $C_{\ken}$ and $C_{\ken}^-$ with
   positive $\mathbb{P}$ probability.
Let $p_j^*=p_j$ if $p_j=0$, and $p_j^*=1$ if $p_j>0$.
 For $k=0,\dots,M\!-\!1$ we denote the expectation matrices
 generated by the vector $(p_0^*,\dots,p_{M-1}^*)$ by $\mathcal{M}^*(k)$.
 Note that if all the level $n$ squares
$Q_{\Ien,\Jen}$ for which $p_{\Ien}>0$ and $p_{\Jen}>0$ are
selected (which happens with positive probability) then
$Z^L(\ken)$ ($Z^R(\ken)$) is the sum of the two elements in the
first column (second column) of $\mathcal{M}^*(\ken)$
respectively. Now note that (with all inequalities componentwise)
\begin{equation}\label{count}
\big(1\; 1\big)\mathcal{M}^*(\ken)\ge \big(1\;
1\big)\mathcal{M}(\ken)\ge \big(\gamma^n\;\gamma^n\big),
\end{equation}
 since the fact that
the sum of the two elements of both of the columns of
$\mathcal{M}(k)$ is greater than $\gamma$ for every $k$ implies
that the sum of the two elements of both of the columns of
$\mathcal{M}(\ken)$ is larger than $\gamma^n$. The claim follows
now directly from (\ref{count}).
\end{proof}

\section{The proof of Theorem 1}

First we show that we can start the $\Delta$-pair process in $R$
at level 2.

\begin{lemma}\label{48} Let $p_{\Delta}$ be the probability that
 $C_{00}\cap \Lambda^2$ contains a level 2 $\Delta$-pair.  If $\gamma>1$
 then $p_{\Delta}>0$.
\end{lemma}

\begin{proof} We know that $\gamma_{M-1}=p_0p_{M-1}+p_1p_0+\cdots
p_{M-2}p_{M-3}+p_{M-1}p_{M-2}>1$. This means that at least one of
the last $M-1$ terms is not zero. So, there is an $0\leq i\leq
M-2$ such that both $p_i>0$ and $p_{i+1}>0$. Then $Q_{i(i+1),ii}$
is selected with probability $p_i^3p_{i+1}$ and
$Q_{(i+1)(i+1),(i+1)(i+1)}$ is selected with probability
$p_{i+1}^{4}$. Using that $L_{i(i+1),ii}=Q_{i(i+1),ii}\cap C_{00}$
and $R_{(i+1)(i+1),(i+1)(i+1)}=Q_{(i+1)(i+1),(i+1)(i+1)}\cap
C_{00}$ we obtain that with probability $p_i^3p_{i+1}p_{i+1}^{4}$
we select the $\Delta $-pair
$(L_{i(i+1),ii},R_{(i+1)(i+1),(i+1)(i+1)})$ in $C_{00}$. Thus
$p_{\Delta }\geq p_i^3p_{i+1}^5>0$.
\end{proof}

It follows from the self-similarity of the construction that the
following fact is true.

\begin{fact}\label{full}
Let $(L^n,R^n)$ be a $n$-th level $\Delta$-pair in some column
$C_{\ken}$. Consider the following
 conditional probability:
 $$\mathbb{P}(\mbox{Proj}_{90^\circ}((L^n\cup R^n)\cap\Lambda)=
 J_{\ken}|\; L^n\subset  \Lambda^n,R^n\subset  \Lambda^n).$$
Then this probability, denoted $p_J$, is independent of $n$,
$\ken$, and the choice of $L^n$ and $R^n$.
\end{fact}

\begin{proposition} \label{23}
Assume that  $\gamma>1$. Then $p_J>0$.
\end{proposition}

In the sequel we will denote
 $$N(\ken ):=\min\left\{Z^L(00\ken  ),Z^R(00\ken )\right\}.$$
Note that $N(\ken )$ counts the  number of level $n+2$
$\Delta$-pairs in subcolumns of $C_{00}$, and that by Lemma~\ref{48}
we know that we start in $C_{00}$ with a level 2  $\Delta$-pair with
positive probability.


 Lemma~\ref{pos} below will directly imply
Proposition ~\ref{23}. In Lemma~\ref{pos} we apply the large
deviation theorem in the same way as in Falconer and Grimmett
(\cite{FG}). Unfortunately in our case the appropriate random
variables are not pairwise independent. To handle this problem we
first prove a lemma which implies that  $N(\ken )$ level $n+2$
$L$-triangles and  $N(\ken )$ level $n+2$ $R$-triangles in
column $C_{00\ken}$ can be paired into $\Delta$-pairs such that
these $N(\ken)$ $\Delta$-pairs can be divided into three groups
(of approximately the same cardinality) with the following
property: any two triangles  (left or right) from any two pairs
from the same group are disjoint. This will provide the required
independence. For each $\ken$ we consider \emph{all} the left and
right triangles in column $C_{00\ken}$. Let $K=K(\ken)$ be their
cardinality. We can naturally label these $K$ triangles with
$\left\{1,\dots ,K\right\}$ in the order in which they appear in
the column, starting at the bottom. Then the odd numbers
correspond to the level $n+2$ $R$-triangles of $C_{00\ken}\cap
\Lambda^{n+2}$ and the even numbers to the $L$-triangles. The
assumption that an $L$ and $R$ triangle form a $\Delta $-pair is
equivalent to the assumption that the corresponding even and odd
numbers are not consecutive. Through this identification the
following combinatorial lemma ensures the division into three
groups announced above.


\begin{lemma}\label{58}
We are given $N$ distinct odd numbers $o_1, \dots, o_N$ and $N$
distinct even numbers $e_1, \dots , e_N$. Then we can couple the
odd numbers with the even numbers and we can color the $N$ couples
with three colors (say $\mathtt{r},\mathtt{g}$ and $\mathtt{b}$)
such that no two numbers in pairs of the same color are adjacent
and all colors are used for at least $\lfloor N/3\rfloor $ pairs.
That is, there exists a permutation $\pi $ of $\left\{1,\dots
,N\right\}$ such that we can color the pairs
$$
(e_1,o_{\pi (1)}),\dots ,(e_N,o_{\pi (N)})
$$
with the three colors such that with each color we painted at least
$\lfloor N/3\rfloor$ pairs and
 for any (also if $\ell=k$)
$(e_k,o_{\pi (k)})$ and $(e_{\ell},o_{\pi (\ell)})$ having the
same color it is true that:
\begin{equation*}
|e_{\ell}-o_{\pi (k)}|>1.
\end{equation*}
\end{lemma}

The proof of this three color lemma will be given in the appendix.


The following key lemma, and its proof, are very similar to the
main result on orthogonal projections of random Cantor sets  in
\cite{FG}.

\begin{lemma}\label{pos}
Assume that $\gamma>1$. Then
\begin{equation*}
  \mathbb{P}\left(N(\ken )>0 \;\forall \ken \in \Madn  \;\mbox{\rm for all }\, n \right)>0
\end{equation*}
holds.
\end{lemma}
\begin{proof}
Using that $\mathbb{E}Z^{V}(k)\geq \gamma$ for $V\in
\{L,R\}$
 and $k\in \Mad$, it
follows from Large Deviation Theory that we can choose an $\eta'$
with $1<\eta' <\min\{2,\gamma\}$ and $0<\delta <1$ such that
\begin{equation}\label{49}
  \mathbb{P}\left(Z_1^V(k)+\cdots +Z_q^V(k)<q\eta' \right)\leq \delta^q
\end{equation}
for all $q\geq 1$, whenever $Z_{1}^{V}(k),Z_{2}^{V}(k)\dots $ are
independent random variables with the same distribution as $Z^V(k)$.
Fix an $1<\eta <\eta '$ and choose $n_0$ such that for all $n\geq n_0$
\begin{equation}\label{67}
\eta \cdot \left(\left\lfloor\frac{\eta ^n}{3}\right\rfloor+1\right)<%
\eta '\cdot \left\lfloor\frac{\eta ^n}{3}\right\rfloor.
\end{equation}

Let $$A_n:=\left\{N(\ken )\geq \eta ^n: \forall \ken \in
\left\{0,\dots, M-1\right\}^n\right\}.$$ It follows from Lemma
\ref{37} and Lemma \ref{48} that for all $n\geq 1$ we have
\begin{equation}\label{50}
  \mathbb{P}(A_n)>0.
\end{equation}
To continue the proof we have to get rid of possible dependence
between $\Delta$-pairs. Fix an arbitrary  $\ken$ and  $k$. Let
$$N:=3\cdot \left\lfloor\frac{N(\ken)}{3} \right\rfloor.$$
Using Lemma~\ref{58} we can we can label and then pair the level
$n+2$ left and right triangles of $C_{00\ken}\cap \Lambda^{n+2}$
into $N$ $\Delta$-pairs $(L_1,R_1),\dots ,(L_N,R_N)$ such that for
every $i=0,1,2$ we have that all the triangles
\begin{equation}\label{57}
L_{iN/3+1},R_{iN/3+1},\dots ,L_{(i+1)N/3},R_{(i+1)N/3} \mbox{ are disjoint.}
\end{equation}
For every $i=0,1,2$ and $1\leq j \leq N/3$ we denote the
$\Delta$-pairs $D_{j}^{i}:=(L_{iN/3+j},R_{iN/3+j})$ and we write
$\widetilde{Z}^V_{iN/3+j}(k)$ for the number of level $n+3$ $V$
triangles in $C_{00\ken k}\cap D_{j}^{i}$. Note that for every
$i$ \ it follows from (\ref{57}) that the $N/3$ random
variables
\begin{equation*}\label{56}
  \widetilde{Z}^V_{iN/3+1}(k),\dots , \widetilde{Z}^V_{(i+1)N/3}(k)
\end{equation*}
are independent and each of them has the same distribution as
$Z^V(k)$. Now we define
$$
S_i^V(00\ken k):=\widetilde{Z}^V_{iN/3+1}(k)+\dots +
\widetilde{Z}^V_{(i+1)N/3}(k).
$$

So, for any $\ken$, $k$ and $V\in \left\{L,R\right\}$ we have
$$
Z^V(00\ken k)\geq \sum\limits_{i=0}^{2}S_i^V(00\ken k).
$$
This directly implies that
\begin{multline*}
 \mathbb{P}\left(Z^V(00\ken k)<\eta^{n+1}|A_n\right)  \leq
\sum\limits_{i=0}^{2} \mathbb{P}\left(S_i^V(00\ken k)<\eta \left(
\left\lfloor\frac{\eta^{n}}{3}\right\rfloor+1\right)\Big|\,A_n\right)\\
  \leq \sum\limits_{i=0}^{2}
\mathbb{P}%
\left(\widetilde{Z}^V_{iN/3+1}(k)+\dots +
\widetilde{Z}^V_{(i+1)N/3}(k)<\eta
     \left(\left\lfloor\frac{\eta^{n}}{3}\right\rfloor+1\right)
    \Big|\,A_n\right),
\end{multline*}
and thus, using that $N\ge \eta^n$ on $A_n$, that the
$Z_{i}^{V}(k)$ are independent of $A_n$, and using (\ref{67}) and
(\ref{49}) we obtain
\begin{eqnarray*}
  \mathbb{P}\left(Z^V(00\ken k)<\eta^{n+1}|A_n\right)
   &\leq &  3\cdot \mathbb{P}\left(
Z_{1}^{V}(k)+\cdots +Z_{\left\lfloor\eta
^n/3\right\rfloor}^{V}(k)<\eta' \left\lfloor\frac{\eta^{n}}{3}\right\rfloor \right)\\
 &\leq &3\delta ^{\lfloor\eta ^n/3\rfloor}.
\end{eqnarray*}

 So
\begin{eqnarray*}
\nonumber  \mathbb{P}\left(A_{n+1}^{c}|A_n\right)%
   &=& %
 \mathbb{P}\Bigg(\bigcup \limits_{V\in
\left\{L,R\right\}} \bigcup \limits_{\ken }\bigcup _k
Z^V(00\ken k)<\eta ^{n+1} \Big|\,A_n\Bigg) \\
  &\leq & \sum\limits_{V\in \left\{L,R\right\}}%
\sum\limits_{\ken }\sum\limits_{k} %
 \mathbb{P}\left(Z^V(00\ken k)<\eta ^{n+1}\Big|\,A_n\right) \\
\nonumber   &\leq & 6\cdot M^n\cdot M\cdot \delta ^{\lfloor\eta
^n/3\rfloor}.
\end{eqnarray*}
Using this and the fact that for any $r\leq n$, we have
$\mathbb{P}\left(A_{n+1}^{c}|A_r\cap \dots \cap
A_n\right)=\mathbb{P}\left(A_{n+1}^c|A_n\right)$ we obtain
that
\begin{equation*}\label{29}
  \mathbb{P}\left(A_{n+1}^c|A_r\cap \dots \cap A_n\right)
  \leq 6M^{n+1}\delta ^{\lfloor\eta ^n/3\rfloor},
\end{equation*}
holds for all $r\leq n$. Therefore for all $r\leq n$,
\begin{eqnarray*}\label{30}
  \mathbb{P}\left(A_{n+1}\cap\dots\cap A_r\right)&\geq&
  \left(1-6M^{n+1}\delta ^{\lfloor\eta ^n/3\rfloor}\right)\mathbb{P}
  \left(A_{n}\cap\dots\cap A_r\right)\\
  &\geq& \mathbb{P}(A_r)
  \prod_{k=r}^{n}\left(1-6M^{k+1}\delta ^{\lfloor\eta ^k/3\rfloor}\right).
\end{eqnarray*}
Choose $r\geq n_0$ such that $\prod_{n=r}^{\infty
}\left(1-6M^{n+1}\delta ^{\lfloor\eta ^n/3\rfloor}\right)>0$. Using
(\ref{50}) this implies that $c:=\mathbb{P}(A_r)\prod_{n=r}^{\infty
}\left(1-6M^{n+1}\delta ^{\lfloor\eta ^n/3\rfloor}\right)>0$. Thus
\begin{equation*}\label{31}
  \mathbb{P}\left(A_n\mbox{ holds for all }n\geq
  r\right)\geq c.
\end{equation*}
This immediately implies the statement of the Lemma.
\end{proof}

\begin{corollary}\label{inter}
Let $\gamma>1$. For every $n$ and $(\Ien,\Jen)$ the conditional
probability of the event that the projection of $\Lambda \cap
Q_{\Ien,\Jen}$ to $J$ contains an interval given that
$Q_{\Ien,\Jen}\subset \Lambda^n$
 is at least $p_{\Delta}\cdot p_J>0$.
\end{corollary}

\begin{proof}
By symmetry it suffices to prove this for squares such that
$R_{\Ien,\Jen}$ is contained in $R$.
  Let $C_{\ken}$ be the column that contains $R_{\Ien,\Jen}$. Then
$C_{\ken 00}$ will contain a $\Delta$-pair of level $n+2$ with
probability not less than $p_{\Delta}$ (see Lemma \ref{48}). Then by
Fact~\ref{full}  and Proposition~\ref{23} the probability that the
projection of this $\Delta$-pair intersected with $\Lambda$ is the
interval $J_{\ken 00}$ is at least $p_J$.
\end{proof}

From here we can finish the proof of our theorem  as in the proof
of Theorem 1 in  \cite{FG}. However, because there is a lot of
dependence between the squares in $\Lambda^{n}$, our proof of part
(a) is slightly more involved.


\begin{proof}[Proof of Theorem~\ref{12} (a)]
Here we assume that $\gamma >1$. We call two squares
$Q_{\Ien,\Jen}$ and $Q_{\Ienn,\Jenn}$ \emph{unaligned} if both
$\Ien\neq \Ienn$ and $\Jen\neq \Jenn$. For every $n$ let $K(n)$ be
the maximal number of pairwise unaligned  squares of $\Lambda^n$.
Then, by maximality of $K(n)$, we can cover $\Lambda $ with $K(n)\cdot 2M^n$ squares of
side $M^{-n}$. So, conditioned on $\Lambda\ne\emptyset$  the Hausdorff dimension of
$\Lambda $ almost surely satisfies
$$
1<\dim_{\rm H}(\Lambda) \le \lim\limits_{n\to\infty}\frac{\log
\left(K(n)\cdot 2M^{n}\right)}{\log M^n}.
$$
Here $\dim_{\rm H}(\Lambda)> 1$ follows from
the hypothesis $\gamma_k>1$ for all $k$ (cf. condition
(\ref{11}), which is equivalent to $\gamma_0+\cdots+
\gamma_{M-1}>M$).

We obtained that
\begin{equation}\label{107}
\{\Lambda \ne\emptyset\} \subset
\{\lim\limits_{n\to\infty}K(n)=\infty\}.
\end{equation}
For every $n$ we fix a system $\{Q^n_1, Q^n_2, \dots,Q^n_{K(n)}\}$ of pairwise unaligned $n$-squares
 contained in $\Lambda^n$ which has cardinality $K(n)$. Let
$$
\mathcal{C}^n_{s}:=%
\{\text{int}( \mbox{Proj}_{90^\circ}(Q^n_s\cap \Lambda))=\emptyset\},
$$
and
$$
\mathcal{C}:=%
\{\text{int}( \mbox{Proj}_{90^\circ}( \Lambda))=\emptyset \},
$$
be the events that the unaligned squares $Q^n_s\cap \Lambda$ for $s=1,\dots,K(n)$, and
the total set $\Lambda$ do not have an interval in their projection.
Our goal is to prove that
\begin{equation*}
\mathbb{P}\left(\mathcal{C}|\,\Lambda \ne\emptyset\right)=0.
\end{equation*}
 According to Corollary \ref{inter} it holds for  all $s$ that
\begin{equation*}
\mathbb{P}\left( \mathcal{C}^n_{s}\right)<1-p_{\Delta}p_J=:t<1.
\end{equation*}
By the definition it is clear that for every $n,N$ we have
\begin{eqnarray*}
\mathbb{P}(\mathcal{C}|\,\Lambda \ne\emptyset)&\le&
        \mathbb{P}(K(n)<N|\,\Lambda \ne\emptyset)+
                \mathbb{P}(\mathcal{C}\cap \{K(n)\ge N\}|\,\Lambda \ne\emptyset)\\
    &\le& \mathbb{P}(K(n)<N|\,\Lambda \ne\emptyset)+
       \mathbb{P}(\mathcal{C}^n_{1}\cap\dots\cap\mathcal{C}^n_{N}|\,\Lambda\ne\emptyset)\\
    &\le& \mathbb{P}(K(n)<N|\,\Lambda
    \ne\emptyset)+\frac{t^N}{\mathbb{P}(\Lambda\ne\emptyset)},
\end{eqnarray*}
where we use that the branching processes  which determine $\Lambda$ in each
 of the unaligned squares  run independently.
Letting first $n\to\infty $, and then $N\to\infty $ we obtain from (\ref{107}) and $t<1$ that
$\mathbb{P}(\mathcal{C}|\,\Lambda \ne\emptyset)=0$ which completes our proof.
\end{proof}

\begin{proof}[Proof of Theorem~\ref{12} (b)]
If there is $k$ such that both $\gamma_k,\gamma_{k+1}<1$ then
using (\ref{36}) we obtain that both of the column sums of the
matrix $\mathcal{M}(k)$ are less than one. This and (\ref{33})
implies that for every $\Ken$ we have
\begin{equation*}
 \lim\limits_{m\to\infty}\|\mathcal{M}
 (\Ken,\underbrace{k,k,...,k}_{m})\|_1 \to 0
\end{equation*}
Let $Z_m$ be the total number of either left or right triangles of
level $n+m$ in column $C_{\Ken,{
\tiny\underbrace{k,k,...,k}_{m}}}$.
  Then
 $$\mathbb{E}(Z_m)=\|\mathcal{M}
 (\Ken,\underbrace{k,k,...,k}_{m})\|_1$$
and by the Markov inequality, $\mathbb{P}(Z_m>0)\leq
\mathbb{E}(Z_m)$. So with probability one  $$\bigcap
_{m=1}^{\infty }C_{\ken,{ \tiny\underbrace{k,k,...,k}_{m}}}$$ is
\emph{not} contained in the projection of $\Lambda $. Since by
varying $\ken$ we may obtain a dense set of such points, we
conclude that in this case the projection of $\Lambda $ does not
contain any interval with probability one.
\end{proof}

\section{Higher order Cantor sets and eigenvalues}
\label{101rho}

Here we reconsider (cf.~Subsection~\ref{scope}) the question of
the scope of Theorem~\ref{12} by introducing higher order Cantor
sets. We also discuss the connection with the eigenvalues of the
matrices involved in the generation of $F_2-F_1$. This will be
illustrated by two examples.
 First we consider the family of random Cantor sets
parametrised by $\rho$ with $0\le\rho\le1$ given by $M=4$ and
$(p_0,\dots,p_3)=(1,0,1,\rho)$. Clearly this gives
  $$
  \mathcal{M}(0)=\!\begin{bmatrix}
    \rho &  0 \\
    \rho  & 2+\!\rho^2
  \end{bmatrix},
  \mathcal{M}(1)=\!\begin{bmatrix}
  1 &  \rho \\
  1  & \rho
\end{bmatrix},
\mathcal{M}(2)=\!\begin{bmatrix}
  \rho &  1 \\
  \rho  & 1
\end{bmatrix},
\mathcal{M}(3)=\!\begin{bmatrix}
2+\!\rho^2     & \rho  \\
0  & \rho
\end{bmatrix}.
  $$
The cyclic autocorrelations are
$$\gamma_0=2+\!\rho^2,\; \gamma_1=2\rho,\; \gamma_2=2,
\;\gamma_3=2\rho.$$ The Palis conjecture predicts that the
difference set will contain an interval almost surely for all
$\rho>0$. Application of Theorem~\ref{12} gives no conclusion for
$\rho<\frac12$, and that for $\rho>\frac12$ this is indeed the
case. However, it is possible to get more out of the theorem by
considering higher order Cantor sets.

The order 2 Cantor set associated to the set generated by
$(p_0,\dots,p_{M-1})$ is the base $M^2$ Cantor set with vector
$$(p_0^{(2)},\dots,p_{M^2-1}^{(2)})$$
given by $$p_{Mi+j}^{(2)}=p_ip_j\quad {\rm{for}}\quad i,j\in
\Mad.$$ We will denote the objects associated to $p^{(2)}$ all
with a superindex $(2)$, for instance $F^{(2)}$ is the random
$M^2$-adic Cantor set generated by $p^{(2)}$, and $I^{(2)}_{\Ken}$
denotes an  $n$-th level $M^2$-adic interval. The key fact is that
for all $\Ien, \Jen \in \Madn$ one has
$$I^{(2)}_{Mi_1+j_1,\dots,Mi_n+j_n}=I_{i_1j_1\dots i_nj_n}.$$
This implies that  $F^{(2)}$ has the same distribution as
$\bigcap_{n\ge 0}F^{2n}$, which equals the original Cantor set
$F$. We can therefore obtain statements about $F$ by applying
Theorem~\ref{12} to $p^{(2)}$. Note that
$\mathcal{M}^{(2)}(Mi+j)=\mathcal{M}(ij)=\mathcal{M}(i)\mathcal{M}(j)$.
So in our example
$$\mathcal{M}^{(2)}(3)=\mathcal{M}(0)\mathcal{M}(3)=
\begin{bmatrix}
    \rho &  0 \\
    \rho  & 2+\!\rho^2
  \end{bmatrix}
\begin{bmatrix}
2+\!\rho^2     & \rho  \\
0  & \rho
\end{bmatrix}=
\begin{bmatrix}
    2\rho+\!\rho^3 &  \rho^2 \\
    2\rho+\!\rho^3  & \rho^2+\!2\rho+\!\rho^3
\end{bmatrix}.
$$
It follows that $\gamma^{(2)}_4=4\rho+2\rho^3$ and
$\gamma^{(2)}_3=2\rho+2\rho^2+\rho^3$. Clearly
$\gamma^{(2)}_3<\gamma^{(2)}_4$, and the latter is smaller than 1
for all $\rho$ smaller than the real root of $4\rho+2\rho^3=1$,
which is about 0.242. Theorem~\ref{12} now gives that $F_2-F_1$
does almost surely not contain an interval for all $\rho<0.242$.
On the other hand we can also strengthen the conclusion for the
opposite case: a straightforward computation yields that for all
$\rho$, $\gamma^{(2)}=2\rho+2\rho^2$. Hence $F_2-F_1$ will contain
an interval for all $\rho$ larger than  $(\sqrt{3}-1)/2
=0.366\dots$.

Note that for all positive $\rho$ the Perron Frobenius eigenvalues
of all the $\mathcal{M}(k)$ are larger than 1, but that still
$F_2-F_1$ does not contain an interval for a range of values of
$\rho$. However, eigenvalues may be useful to prove the opposite
case: the Perron Frobenius eigenvalue of $\mathcal{M}^{(2)}(3)$ is
equal to
$$\left(\rho^2+\rho/2+2+1/2 \sqrt{4\rho^3+\rho^2+8\rho}\right)\rho,$$
which is smaller than 1 when $\rho<0.3221$. As in the proof of
Theorem~\ref{12} part(b), this can be used to show that a dense
set of points is not in
 $F_2-F_1$.
 Using higher order Cantor sets (up to order 324), and Matlab we obtained that
 the critical point $\rho_c$ where  $F_2-F_1$ changes
 from empty to non empty interior, satisfies  $0.3222<\rho_c<0.3226$.

\begin{figure}[h!]
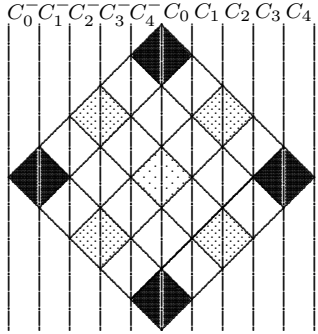

$$
\beginpicture
\setcoordinatesystem units <0.03\textwidth,0.03\textwidth>
\setplotarea x from -1 to 9, y from -5 to 5.6
\setlinear
\linethickness=13pt
\setlinear
\plot  0 0    5 5  /
\plot  0 0   5  -5 /
\plot  4.03  -3.97  8 0  /
\plot  1 -1    6 4  /
\plot  2 -2    7 3  /
\plot  3 -3    8 2  /
\plot  4 -4    9 1  /
\plot  5 -5    10 0  /
\plot  0 -5    0  5 /
\plot  1 -5    1  5 /
\plot  2 -5    2  5 /
\plot  3 -5    3  5 /
\plot  4 -5    4  5 /
\plot  5 -5    5  5 /
\plot  6 -5    6  5 /
\plot  7 -5    7  5 /
\plot  8 -5    8  5 /
\plot  9 -5    9  5 /
\plot  10 -5    10 5 /
\plot  1  1    6  -4  /
\plot  2  2    7  -3  /
\plot  3  3    8  -2  /
\plot  4  4    9  -1  /
\plot  5  5    10  0  /
\setshadegrid span <0.01cm> 
\vshade 0 0 0 1 -1 1 2 0 0 /
\setshadegrid span <0.075cm> 
\vshade 4 0 0 5 -1 1 6 0 0 /
\setshadegrid span <0.050cm> 
\vshade 2 2 2 3 1 3 4 2 2 /
\setshadegrid span <0.01cm> 
\vshade 4 4 4 5 3 5 6 4 4 /
\setshadegrid span <0.05cm> 
\vshade 2 -2 -2 3 -3 -1 4 -2 -2 /
\setshadegrid span <0.01cm> 
\vshade 4 -4 -4 5 -5 -3 6 -4 -4 /
\setshadegrid span <0.01cm> 
\vshade 8 0 0 9 -1 1 10 0 0  /
\setshadegrid span <0.05cm> 
\vshade 6 2 2 7 1 3 8 2 2 /
\setshadegrid span <0.05cm> 
\vshade 6 -2 -2 7 -3 -1 8 -2 -2 /
{\tiny
\put {$C_{0}$} at 5.5 5.4
\put {$C_{1}$} at 6.5 5.4
\put {$C_{2}$} at 7.5 5.4
\put {$C_{3}$} at 8.5 5.4
\put {$C_{4}$} at 9.5 5.4
\put {$C_{0}^-$} at 0.5 5.4
\put {$C_{1}^-$} at 1.5 5.4
\put {$C_{2}^-$} at 2.5 5.4
\put {$C_{3}^-$} at 3.5 5.4
\put {$C_{4}^-$} at 4.5 5.4
} 
\endpicture
$$
\caption{The square $Q$ for the  $(1,0,\rho,0,1)$ Cantor set; the
shades indicate the probabilities with which the $Q_{i,j}$ do
occur.} \label{fig:FG}
\end{figure}

We consider a second parametrized family which has a very different
behaviour. Let $M=5$, and $(p_0,\dots,p_4)=(1,0,\rho,0,1)$ for
$0\le\rho\le1$. (See Figure \ref{fig:FG}.)
 One finds
  $$
  \mathcal{M}(0)=\!\begin{bmatrix}
   1 &  0 \\
   0  & 2+\!\rho^2
  \end{bmatrix},
  \mathcal{M}(1)=\!\begin{bmatrix}
  0 &  1 \\
  2\rho  & 0
\end{bmatrix},
\mathcal{M}(2)=\!\begin{bmatrix}
  2\rho &  0 \\
  0 &  2\rho
\end{bmatrix},
  $$
and $\mathcal{M}(3), \mathcal{M}(4)$ are obtained from
$\mathcal{M}(1)$, respectively $\mathcal{M}(0)$ by interchanging
$R$ and $L$. Since $\gamma^{(n)}_1=1$ for all $n$,
Theorem~\ref{12} is not applicable, even if one considers higher
order Cantor sets. However, since the 5 matrices are either
diagonal or anti-diagonal, it is not hard to prove that the
Perron-Frobenius eigenvalues $\lambda(\Ken)$ of the matrices
$\mathcal{M}(\Ken)$  satisfy
$$\lambda(\Ken)\ge (\sqrt{2\rho})^n.$$
This seems to suggest that the critical $\rho$ for this family is
equal to 1/2. Surprisingly, we here have $\rho_c=1$. It follows
from Theorem~\ref{theorem2} that $F_2-F_1$ contains an interval
when $\rho=1$. To see that $F_2-F_1$ has empty interior for
$\rho<1$, consider column $C_{\Ien 44\cdots42}$ of level $n+m+1$
for each $\Ien$ such that $C_{\Ien}$ has only $R$-triangles (these
are in fact all the columns where $i_1,\dots ,i_n$ are all even
numbers). Then if there are $K$ $R$-triangles in $C_{\Ien}$ this
will be also true for all columns $C_{\Ien 44\cdots4}$ of level
$n+k$, where $k=1,2,\dots,m$, and moreover, column $C_{\Ien
44\cdots42}$ will be empty with probability
$$\big[(1-\rho)^2\big]^K$$ for all $m$. It follows as at the end of
the proof of Theorem~\ref{12} b) that there is a dense set of points
in the complement of the projection of $\Lambda$.

An alternative would be to adapt the  proof (to our setting which
has much more dependence) of the main result on fractal
percolation of \cite{FG}, or rather its supplement from
\cite{FGcorr}. The crucial observation here is that the first
level columns $C_k^-$ and $C_k$  split in pairs
$(C_{0}^-,C_{1}^-),\dots,(C_{3},C_{4}) $, that do not interact
with each other, since the $Q_{i,j}$ only occur with positive
probability when $i-\!j$ is even. Redefining each pair of columns
as a single new column, and tilting the $Q_{i,j}$'s, the question
of empty interior could then be resolved by the (extended) results
of \cite{FG} and \cite{FGcorr}.

\section{The deterministic case}
\label{determ} In the deterministic case each $p_i$ is either 0 or
1, and the matrices $\mathcal{M}(k)$ simply count the number of
level 1 $L$-triangles and $R$-triangles in the columns $C_k^-$ and
$C_k$. The crux to the solution in this case is that we can reduce
the problem to a finite problem by observing that to have a
non-empty projection in a certain column we only have to know
whether there is \emph{at least one} $L$- or $R$-triangle in that
column. We relax the $\Delta$-pair condition: now a pair of
adjacent $R$- and $L$-triangles is also allowed since independence
is no longer an issue.


For a non-negative integer matrix $A$, let its \emph{reduction}
$A^\triangledown$ be defined by $a_{ij}^\triangledown=0$
 if $a_{i,j}=0$, and $a_{ij}^\triangledown=1$ if
$a_{i,j}\ge 1$. Note that the reduction of the product of two
matrices equals the product of their reductions. It follows that the
reduction of
$\mathcal{M}(k_1)^\triangledown\cdots\mathcal{M}(k_n)^\triangledown$
describes the presence or absence of $n$-th order $L$- and
$R$-triangles in columns $\Ken$ of order $n$. Let $\mathcal{T}$ be
the set of $2\times 2$ matrices with entries 0 and 1. For
convenience we denote these matrices by their natural binary coding:
$$\begin{bmatrix}
    a &  b \\
    c  & d
  \end{bmatrix}
=T_j\quad \Leftrightarrow \quad j=a+2b+4c+8d.$$ Define the map
$G:2^\mathcal{T}\rightarrow  2^\mathcal{T}$ by
$G(\emptyset)=\emptyset$, and for $\mathcal{C}\ne\emptyset$
$$G(\mathcal{C})=\{(TT')^\triangledown: T\in \mathcal{C}, T'\in\mathcal{C}\}.$$
Then there is an empty column of order $n$ in $\Lambda^n$ if and
only if
$$T_0=
\begin{bmatrix}
  0 &  0 \\
  0  & 0
\end{bmatrix}
\in
G^n(\{\mathcal{M}(0)^\triangledown,\cdots,\mathcal{M}(M\!-\!1)^\triangledown\}),
$$
where $G^n$ is the $n$-th iterate of $G$. Since $G$ is acting on a
finite set, the orbit of any point  becomes eventually periodic.
We call that periodic sequence
 an \emph{attractor}, denoted by  $\mathcal{A}$. Examples are fixed
 points of $G$, as e.g.,  $\mathcal{A}=\{T_0\}$ and
 $\mathcal{A}=\{T_6,T_9\}$. Assisted by the computer we can show
 that actually \emph{all} attractors are fixed points (the proof
 below can be adapted to a proof which does not explicitly use
 this result).

\begin{theorem}\label{theorem2}
Let the Cantor set $F$ be generated by a 0-1-vector
$(p_0,\dots,p_{M-1})$.
 Then   $F_2-F_1$  does not contain any intervals if and only if
 $T_0\in \mathcal{A}$, where $\mathcal{A}$ is the fixed point of
 the map $G$ starting from
 $\{\mathcal{M}(0)^\triangledown,\cdots,\mathcal{M}(M\!-\!1)^\triangledown\}$.
\end{theorem}

\begin{proof}
$\Leftarrow)$ \quad If $T_0\in \mathcal{A}$ then an empty column
has to occur in some column of order $n$, where $n\le 2^{16}$
(actually a computer analysis shows that $n\le 3$). The proof that
$F_2-F_1$  does not contain any intervals, is then finished as the
proof of Theorem~\ref{12}, part b).

$\Rightarrow)$ \quad Suppose  $T_0\notin \mathcal{A}$. We split
into two cases.

\smallskip

\textit{Case 1. For some $n\ge 1$ a $\Delta$-pair of order $n$ occurs.}

\smallskip

Suppose that this happens in column $C_{\ken}$ or in $C^-_{\ken}$.
For arbitrary $m$ and ${l_1\dots l_m}$ the fact that $T_0$ does
not occur in $\mathcal{A}$  implies that $\mathcal{M}(l_1\dots
l_m)\ne T_0$, and hence that all subcolumns $C_{\ken\lem}$ will
contain at least one order $n+m$ triangle for all $m$, and so the
complete interval $J_{\ken}$ respectively $J^-_{\ken}$ will lie in
the projection of $\Lambda$.

\smallskip

\textit{Case 2. A $\Delta$-pair never occurs.}

\smallskip

Then $\mathcal{A}$ can not contain a matrix with a row of two 1's. This means that
$$\mathcal{A}\subset \{T_1,T_2,T_4,T_5,T_6,T_8,T_9,T_{10}\}.$$
But since $T_2^2=T_4^2=T_0$, these two matrices can also not
occur, and hence
$$\mathcal{A}\subset \{T_1,T_5,T_6,T_8,T_9,T_{10}\}.$$
Now suppose that $T_1\in \mathcal{A}$.
Then, since $T_1T_8=T_0$, $T_1T_6=T_4$ and $T_1T_{10}=T_0$, it follows (using again that
$T_4^2=T_0$) that
$$\mathcal{A}\subset \{T_1,T_5,T_9\}, \qquad \mbox{given that}\;  T_1\in \mathcal{A}.$$
Now note that all three matrices  $T_1,T_5$, and $T_9$ have a 0 in
the $LR$ position. This implies that for a certain $n$ (actually
it is not hard to show that one can take $n=1$) the $n$-th order
Cantor set $\Lambda^n$ has the property that there are no
$R$-triangles in its intersection with the $L$ triangle. This
happens only if at most one $p_i^{(n)}\ne 0$, which contradicts
(\ref{nonem}). Conclusion:  $T_1\notin \mathcal{A}$. Analogously
(replacing $L$ by $R$), it follows that $T_8\notin \mathcal{A}$.
So we find that necessarily
$$\mathcal{A}\subset \{T_5,T_6,T_9,T_{10}\}.$$
But the matrices $T_5,T_6,T_9$, and $T_{10}$ each have at least
one 1 in each row. It follows that for \emph{all} $n$ \emph{all}
columns of order $n$ contain at least one triangle, i.e., that
$\mbox{Proj}_{90^\circ}(\Lambda)$ is the whole interval
 $[-\frac12\sqrt{2},\frac12\sqrt{2}]$.
\end{proof}

\section{Appendix: proof of Lemma~\ref{58}}

\begin{proof} Let
$S_0:=\left\{o_1, \dots, o_N,e_1, \dots, e_N\right\}$. We say that
$S_0$ is a \emph{3C set} if the assertion of the Lemma holds for
$S_0$. First we prove that:

\begin{center}
(*) If $S_0$ consists of $2N$ consecutive numbers
$S_0:=\left\{u_1,\dots ,u_{2N}\right\}$\\ then $S_0$ is a 3C set.
\end{center}
To see this we write $N=3p+r$ where $r=0,1$ or $2$. Then
we couple and color the first $6p$ numbers of $S_0$ as follows:
\begin{eqnarray*}
  (u_{3k+1},u_{3k+4}) &=& (\mathtt{r}, \mathtt{r}) \\
 (u_{3k+2},u_{3k+5}) &=& (\mathtt{g}, \mathtt{g})\\
 (u_{3k+3},u_{3k+6}) &=& (\mathtt{b}, \mathtt{b}),
\end{eqnarray*}
for $0\leq k\leq p-1$. Since $\lfloor N/3\rfloor=\lfloor
p+r/3\rfloor=p $  we have verified the assertion of (*) in this way
(without having actually colored the last $2r$ numbers,
an option which we will leave open till the end of the proof). \\
 A subset $B\subset \mathbb{N}$ is \emph{connected} if $n_1,n_2\in B$
and $n_1\leq k\leq n_2, k\in \mathbb{N}$ implies that $k\in B$. We
say that a subset $I\subset S_0$  is an \emph{interval} of $S_0$ if $I$
is a maximal connected subset of $ S_0$ (it is allowed that $I$ consists of one
element). In particular if
 $J_1$ and $J_2\neq J_1$ are intervals of $S_0$ then there exists an
 $\ell\not\in  S_0$
such that $\ell$ separates $J_1$ and $J_2$. Let
$\mathcal{I}^{eo}_0$ be the family of the intervals of $S_0$ for
which the left endpoint is an even number and the right endpoint
is an odd number. Analogously we define the family of intervals
$\mathcal{I}^{ee}_0$, $\mathcal{I}^{oe}_0$ and
$\mathcal{I}^{oo}_0$. Let $\mathcal{I}_0$ be the family of all of
these intervals. So,
$$
S_0=\bigcup _{I\in \mathcal{I}^{ee}_0}I\cup%
\bigcup _{I\in \mathcal{I}^{eo}_0}I\cup%
 \bigcup _{I\in\mathcal{I}^{oe}_0}I\cup%
 \bigcup _{I\in \mathcal{I}^{oo}_0}I=\bigcup \limits_{I\in \mathcal{I}_0}I.
$$
Let us define a ``gluing and shifting" operation $\Phi $ on
$\mathcal{I}_0$ as follows: if there exist two intervals
$J_i=[k_i,\ell_i]\in \mathcal{I}_0$, $i=1,2$ such that $\ell_1+k_2=1
\mod2$ then we select the two left most intervals with this property
and  we form the interval
$$J:=(J_1+n)\cup (J_2+\ell_1-k_2+1+n),$$
where $n\in \mathbb{N}$ is the smallest number such that $J$ is
separated by a distance of at least $2$ from any intervals of
$\mathcal{I}_0\setminus\left\{J_1,J_2\right\}$. In this case we define
$$\mathcal{I}_1:=\Phi (\mathcal{I}_0):=\left\{J\right\}\cup
\mathcal{I}\setminus{\left\{J_1,J_2\right\}} .$$  If there are no
such $J_1,J_2$ then let $\mathcal{I}_1:=\Phi
(\mathcal{I}_0):=\mathcal{I}_0$. By induction we define
$\mathcal{I}_k$ for every $k$. We obtain also by induction that
the set
$$
S_k:=\bigcup _{I\in \mathcal{I}_k}I
$$
consists  of $N$ odd numbers and $N$ even numbers and if $S_{k}$
is a 3C set then $S_{k-1}$ is also a 3C set for
all $k\geq 1$. We claim that
\begin{equation}\label{60}
\mbox{ if }%
\#(\mathcal{I}_k)\geq 3%
\mbox{ then } \#(\mathcal{I}_{k+1})<\#(\mathcal{I}_k).
\end{equation}
We argue by contradiction. If
$\#(\mathcal{I}_{k+1})=\#(\mathcal{I}_k)$ then
$\mathcal{I}_{k+1}=\mathcal{I}_k$. Observe that
\begin{equation}\label{61}
\mathcal{I}_{k+1}=\mathcal{I}_k \mbox{ implies that }%
\mathcal{I}_{k}^{ee}=\mathcal{I}_{k}^{oo}=\emptyset.
\end{equation}
Namely, the cardinality of
$\mathcal{I}_{k}^{ee},\mathcal{I}_{k}^{oo}$  is the same since we
have $N$ odd numbers and $N$ even numbers in $S_k$ and if this
cardinality is not $0$ then we can form an interval $J$ like above
by selecting a $J_1$ from $\mathcal{I}_{k}^{oo}$ and a $J_2$ from
$\mathcal{I}_{k}^{ee}$. Further, if either
$\mathcal{I}_k^{eo}\ne\emptyset $ or
$\mathcal{I}_k^{oe}\ne\emptyset $ then their cardinality is at
most one.  This is so because otherwise choosing $J_1\ne J_2$ from
the same family, we could form $J$ like above. In this way we have
verified (\ref{60}). Let $k_0$ be then smallest number for which
$\mathcal{I}_{k_0}=\mathcal{I}_{k_0+1}$. Then it follows from
(\ref{61}) that
\begin{equation*}
\mathcal{I}_{k_0}=\mathcal{I}_{k_0}^{eo}\cup
\mathcal{I}_{k_0}^{oe}
\end{equation*}
where both of the families on the right hand side consists of at
most one interval. Let us call these intervals $I_1$ and $I_2$
with the possibility that one of them  may be empty. That is
\begin{equation*}
  S_{k_0}=I_1\cup I_2.
\end{equation*}
Since both of the intervals $I_1,I_2$ consist of  an even number of
consecutive numbers, it follows from (*) that they are both 3C
sets. This implies that $S_{k_0}$ is also a 3C set.
The only thing to check is that
\begin{equation}\label{68}
\mbox{ we use all  colors at least }\lfloor N/3\rfloor \mbox{ times }
\end{equation}
since $I_1$ and $I_2$ are separated by a distance of at least $2$.
To see that we can accomplish this, write   $N_i=3p_i+r_i$ where
$N_i$ is the cardinality of $I_i$, and $r_i=0,1$ or $2$. Now if
$r_1+r_2\le 2$, then $N=N_1+N_2=3(p_1+p_2)+r_1+r_2$ and  $\lfloor
N/3\rfloor=p_1+p_2$, so (\ref{68}) is fulfilled.  What remains are
the cases $r_1=1, r_2=2$ (together with  $r_1=2, r_2=1$ which is
very similar), and $r_1=2, r_2=2$, in which cases $\lfloor
N/3\rfloor=p_1+p_2+1$. In the first case we color the last two
numbers  of $I_1$ by  $\mathtt{g}$ and $\mathtt{b}$, and the last
four numbers of $I_2$ by
$\mathtt{r},\mathtt{b},\mathtt{g},\mathtt{r}$. Since the parities
of these numbers are $e,o$, respectively $o,e,o,e$, we see that we
can create the required extra pair for each color. In the case
$r_1=2, r_2=2$ we color the last four numbers of both $I_1$ and
$I_2$ by $\mathtt{r},\mathtt{g},\mathtt{b},\mathtt{r}$, and again
we can create an extra pair for each color.  This proves
(\ref{68}).  As we observed above, the fact that $S_{k_0}$ is a 3C
set implies that $S_0$ is also a 3C set, and the proof is
complete.
\end{proof}


\begin{center}
{\sc Acknowledgement}
\end{center}

The authors would like to say thanks to an anonymous referee for a
number of comments which put us on the right track regarding the
dependence problems.
 The first author is grateful for fruitful discussions with Wenxia Li.
Also the second author would like to say thanks to Michal Rams for
his useful comments.


\end{document}